\documentclass[12pt,a4paper]{article}
\usepackage{epsfig,latexsym,amsfonts,amssymb,amsmath,amscd,
epic,graphics,theorem}
\font\bcyr=wncyb10 
\font\cyri=wncyi10 
\theoremstyle{plain}
\newtheorem{lemma}{Lemma}
\newtheorem{proposition}[lemma]{Proposition}
\newtheorem{remark}[lemma]{Remark}

\newtheorem{theorem}[lemma]{Theorem}

{\theorembodyfont{\rmfamily}  \font\ncsc=cmcsc10
 \font\ntt=cmtt12

\begin{document}
\newcommand{\pperp}{\hbox{$\perp\hskip-6pt\perp$}}
\newcommand{\ssim}{\hbox{$\hskip-2pt\sim$}}
\newcommand{\aleq}{{\ \stackrel{3}{\le}\ }}
\newcommand{\ageq}{{\ \stackrel{3}{\ge}\ }}
\newcommand{\aeq}{{\ \stackrel{3}{=}\ }}
\newcommand{\bleq}{{\ \stackrel{n}{\le}\ }}
\newcommand{\bgeq}{{\ \stackrel{n}{\ge}\ }}
\newcommand{\beq}{{\ \stackrel{n}{=}\ }}
\newcommand{\cleq}{{\ \stackrel{2}{\le}\ }}
\newcommand{\cgeq}{{\ \stackrel{2}{\ge}\ }}
\newcommand{\ceq}{{\ \stackrel{2}{=}\ }}
\newcommand{\N}{{\mathbb N}}
\newcommand{\A}{{\mathbb A}}
\newcommand{\K}{{\mathbb K}}
\newcommand{\Z}{{\mathbb Z}}
\newcommand{\R}{{\mathbb R}}
\newcommand{\C}{{\mathbb C}}
\newcommand{\Q}{{\mathbb Q}}
\newcommand{\PP}{{\mathbb P}}
\newcommand{\Pic}{{\operatorname{Pic}}}\newcommand{\Sym}{{\operatorname{Sym}}}
\newcommand{\oeps}{{\overline\eps}}\newcommand{\idim}{{\operatorname{idim}}}
\newcommand{\oDel}{{\widetilde\Del}}
\newcommand{\real}{{\operatorname{Re}}}\newcommand{\Aut}{{\operatorname{Aut}}}
\newcommand{\conv}{{\operatorname{conv}}}\newcommand{\Id}{{\operatorname{Id}}}
\newcommand{\Span}{{\operatorname{Span}}}
\newcommand{\Ker}{{\operatorname{Ker}}}
\newcommand{\Ann}{{\operatorname{Ann}}}
\newcommand{\Fix}{{\operatorname{Fix}}}
\newcommand{\sign}{{\operatorname{sign}}}
\newcommand{\Tors}{{\operatorname{Tors}}}
\newcommand{\alg}{{\operatorname{alg}}}
\newcommand{\oi}{{\overline i}}
\newcommand{\oj}{{\overline j}}
\newcommand{\ob}{{\overline b}}
\newcommand{\os}{{\overline s}}
\newcommand{\oa}{{\overline a}}
\newcommand{\oy}{{\overline y}}
\newcommand{\ow}{{\overline w}}
\newcommand{\ot}{{\overline t}}
\newcommand{\oz}{{\overline z}}
\newcommand{\eps}{{\varepsilon}}
\newcommand{\proofend}{\hfill$\Box$\bigskip}
\newcommand{\Int}{{\operatorname{Int}}}
\newcommand{\pr}{{\operatorname{pr}}}
\newcommand{\Hom}{{\operatorname{Hom}}}
\newcommand{\rk}{{\operatorname{rk}}}\newcommand{\Ev}{{\operatorname{Ev}}}
\newcommand{\im}{{\operatorname{Im}}}
\newcommand{\sk}{{\operatorname{sk}}}\newcommand{\DP}{{\operatorname{DP}}}
\newcommand{\const}{{\operatorname{const}}}
\newcommand{\Sing}{{\operatorname{Sing}}\hskip0.06cm}
\newcommand{\conj}{{\operatorname{Conj}}}
\newcommand{\Cl}{{\operatorname{Cl}}}
\newcommand{\Crit}{{\operatorname{Crit}}}
\newcommand{\Ch}{{\operatorname{Ch}}}
\newcommand{\discr}{{\operatorname{discr}}}
\newcommand{\Tor}{{\operatorname{Tor}}}
\newcommand{\Conj}{{\operatorname{Conj}}}
\newcommand{\vol}{{\operatorname{vol}}}
\newcommand{\defect}{{\operatorname{def}}}
\newcommand{\codim}{{\operatorname{codim}}}
\newcommand{\ov}{{\overline v}}
\newcommand{\ox}{{\overline{x}}}
\newcommand{\bw}{{\boldsymbol w}}
\newcommand{\bx}{{\boldsymbol x}}
\newcommand{\bd}{{\boldsymbol d}}
\newcommand{\bz}{{\boldsymbol z}}\newcommand{\bp}{{\boldsymbol p}}
\newcommand{\tet}{{\theta}}
\newcommand{\Del}{{\Delta}}
\newcommand{\bet}{{\beta}}
\newcommand{\kap}{{\kappa}}
\newcommand{\del}{{\delta}}
\newcommand{\sig}{{\sigma}}
\newcommand{\alp}{{\alpha}}
\newcommand{\Sig}{{\Sigma}}
\newcommand{\Gam}{{\Gamma}}
\newcommand{\gam}{{\gamma}}
\newcommand{\Lam}{{\Lambda}}
\newcommand{\lam}{{\lambda}}
\newcommand{\SC}{{SC}}
\newcommand{\MC}{{MC}}
\newcommand{\nek}{{,...,}}
\newcommand{\cim}{{c_{\mbox{\rm im}}}}
\newcommand{\clM}{\tilde{M}}
\newcommand{\clV}{\bar{V}}

\title{Welschinger invariants revisited}
\author{Ilia Itenberg
\and Viatcheslav Kharlamov \and Eugenii Shustin}
\date{}
\maketitle
\bigskip

\centerline{\small{\it To the memory of Mikael Passare, remarkable mathematician and beautiful personality}}

\bigskip
\begin{abstract}
We establish the enumerativity of (original and
modified)
Welschinger invariants for every real divisor on any real algebraic
del Pezzo surface and
give an algebro-geometric proof of the invariance of that count both up to variation of the point
constraints on a given surface and variation
of the complex structure of the surface itself.
\end{abstract}

{\hfill {\bcyr
- My govorim s tobo\rm\u{\bcyr i}\bcyr\ na raznyh yazykah, kak vsegda, - otozvalsya Voland,}}

\vskip-5pt

{\hfill {\bcyr - no vewi, o kotoryh my govorim, ot e1togo ne menyayut{}sya.}}

\vskip5pt

{\sc\small \hfill {\cyri M. Bulgakov. Master i Margarita.}}\footnote{"We speak different languages, as usual," responded Woland, "but this does not change the things we speak about. "- {\it M. Bulgakov. The Master and Margarita.}}

\vskip-5pt

\section*{Introduction}

The discovery of Welschinger invariants \cite{W,W1} has revolutionized
real
enumerative geometry. Since then much effort was devoted to the numerical study of Welschinger invariants,
especially in the case of real
del Pezzo surfaces, which allowed one to
prove
long time stated conjectures on existence of real
solutions in corresponding enumerative problems and to observe a new, unexpected phenomena of abundance (see \cite{ABL,HS,IKS,IKS1,IKS2,M});
it also led to introducing certain
modified
Welschinger invariants (see \cite{IKS1}). This development raised several natural questions:
first, for which real
del Pezzo surfaces the Welschinger invariants are
strongly
enumerative
({\it i.e.}, provided by a count, with weights $\pm$ 1,
of real rational curves in a given divisor class, passing through a suitable number of real and complex conjugated points) and, second,
to
what extent such a count is invariant under deformations of the surface.
The enumerative nature of the invariants in the symplectic setting is the key point of \cite{W1}, but it does not imply their enumerative nature
in the algebro-geometric setting because of stronger genericity assumptions. The deformation invariance in the symplectic setting implies the deformation invariance in the
algebro-geometric
setting,
but in \cite{W1}
the symplectic deformation invariance is declared
without proof.
Therefore,
our
principal motivation has been to answer
the question on
algebro-geometric
enumerativity of Welschinger invariants on
real del Pezzo surfaces, and
to prove the deformation invariance in the
algebro-geometric
setting. Our second motivation is
an expectation
that a good understanding of enumeration of real rational curves on real del Pezzo surfaces can help
to extend the results to other types of surfaces and to curves of higher genus (such an expectation is
confirmed now by \cite{IKSrel,ShustinNew}).
The algebro-geometric framework can be also helpful in the study of algorithmic and complexity aspects.

In most of the papers on the subject, the
algebro-geometric
enumerativity of Welschinger invariants on del Pezzo surfaces is considered as known.
Indeed,
it follows from enumerativity of Gromov-Witten invariants for such surfaces, and in the literature on Gromov-Witten invariants the latter enumerativity is considered as known. However, a careful analysis, see Lemma \ref{ln12}, has shown that there is one, and luckily only one,
exception (apparently
not mentioned in the literature): that is the case of del Pezzo surfaces of canonical degree $1$
and $D=-K_\Sig$; for any other pair of a real del Pezzo surface
and a real divisor on it, the Welshinger invariants, original and modified,
are strongly enumerative (in the above exceptional case, the number of solutions is still finite, but certain
solutions may acquire some nontrivial multiplicity).

To prove the deformation invariance, we split the task into two parts.
First, we fix the complex structure and vary the position of the points.
Here, our strategy is close to that of
the original proof of Welschinger in \cite{W1}, but uses
algebro-geometric
tools instead of symplectic ones.
In fact, already some time ago in \cite{IKSa} we have undertaken an attempt to give a purely
algebro-geometric
proof of such an invariance.
However, that proof appears to be incomplete, since one type of local bifurcations
in the set of counted curves was missing;
it shows up for del Pezzo surfaces of canonical degree 1
and $D=-2K_\Sig$ (see Lemma \ref{ln16}(i) below, which states, in particular,
that the closure of the (one-dimensional)  family of rational curves in $|-2K_\Sigma|$ contains non-reduced curves).
To the best of our knowledge, up to now this bifurcation has not been addressed in the literature,
but it is unavoidable even in the symplectic setting (contrary to \cite[Remark 2.12]{W1}).
This step is summed up in Proposition \ref{p1}, which
states the invariance of the Welschinger count under the variation of points for any real divisor on each real del Pezzo surface.

The crucial point of the next step is
the invariance under crossing the walls that correspond to, so-called,  uninodal  del Pezzo surfaces.
Here, our proof is based on a real version of the Abramovich-Bertran-Vakil formula (note that adapting the formula to the symplectic setting one can prove the symplectic deformation invariance following the same lines). In addition, as in the study of the enumerativity,
there
appears
a case not to miss and to investigate separately,
here this is the case of del Pezzo surfaces of canonical degree $1$
and $D=-K_\Sig$.

The paper is organized as follows.
In Section 1 we recall a few basic facts concerning del Pezzo surfaces and their deformations,
introduce Welschinger invariants
in their
modified
version
and formulate the main results. Section 2 develops technical tools needed for the proof of the main results.
There, we study moduli spaces
of stable maps of pointed genus zero curves to del Pezzo
surfaces
and uninodal del Pezzo
surfaces,
describe
generic elements of these moduli spaces and
generic elements of
the
codimension one strata. We show also that Welschinger numbers extend by continuity from the case of immersions
to the case of birational stable maps with arbitrary singularities.
Section 3 is devoted to the proof of the main results.

\vskip5pt

{\bf Acknowledgements}
An essential part of the work was done during the authors' stay at
the Centre Interfacultaire Bernoulli at EPFL in the framework of the program
``Tropical geometry in its complex and symplectic aspects".
The authors are grateful to CIB for the support
and excellent working conditions.
The third author has also been supported from the German-Israeli
Foundation grant no. 1174-197.6/2011 and from the
Hermann-Minkowski-Minerva Center for Geometry at the Tel Aviv
University. We are grateful to G.-M. Greuel and I. Tyomkin for
helpful discussions of certain aspects of deformations of curves.

\section{Definitions and main statements}\label{introduction}

\subsection{Surfaces under consideration}
Over $\C$, a del Pezzo surface is either
$(\PP^1)^2$ or
$\PP^2$ blown up at $0\le k\le 8$ points. Conversely, blowing up
$0\le k\le 8$ points of $\PP^2$ yields a del Pezzo surface if and
only if no $3$ points lie on a straight line, no $6$ lie on a conic, and no $8$ points
lie on a rational cubic
having a singularity at one of these $8$ points.

Del Pezzo
surfaces of degree $d=K^2=9-k\ge 5$ have no moduli. If $d=9$ or
$7\ge d\ge 5$, then there is only one, up to isomorphism, del Pezzo surface of degree $d$
and it can be seen as a blown up $\PP^2$. If $d=8$, then there are 2 isomorphism classes: $(\PP^1)^2$ and $\PP^2$ blown up at a point.
The latter two surfaces are not deformation equivalent.
For
$4\ge d\ge 1$ the moduli space of del Pezzo surfaces of degree
$d=9-k$ is an irreducible $(2k-8)$-dimensional variety.

All
del Pezzo surfaces of given degree $d\ne 8$ are deformation equivalent
to each other, and, for our purpose, it will be more convenient to use, instead of
the moduli spaces, the deformation spaces,
that is, to fix in each deformation class one of the del Pezzo surfaces
(say, a blow up of $\PP^2$ at a certain generic collection of points) and consider the Kodaira-Spencer-Kuranishi space,
{\it i.e.},
the space of all
complex structures on the underlying smooth $4$-manifold
factorized by the action of diffeomorphisms isotopic to identity. Naturally, we awake this
space only when $d\le 4$. We denote it by $\mathcal D_d $.
Del Pezzo surfaces of degree $d$ form
in $\mathcal D_d$ an open dense subset, which we denote by $\mathcal D_d^{DP}$.

The problem of deformations of complex structures on rational surfaces is not obstructed, since $H^2(X,{\mathcal T}_X)=0$ for any smooth rational surface
$X$ (here and further on, we denote by ${\mathcal T}_X$ the tangent sheaf).
In addition, for degree $d\le 4$
del
Pezzo surfaces as well as
for
any generic smooth rational surface $X$ with $K^2_X\le 4$,
we have $H^0(X, {\mathcal T}_X)=0$,
so that at such points the Kodaira-Spencer-Kuranishi space is smooth (but not necessarily Hausdorff).

In fact, the only properties of this space which we use further on  are the following.
We call a surface $\Sig\in{\mathcal D}_d$
{\it uninodal del Pezzo}
if it contains a smooth rational $(-2)$-curve $E_\Sig$, and $-K_\Sig C>0$
for each irreducible curve
$C \ne E_\Sig$
(in particular, $C^2\ge-1$). For $d\le 4$, denote by
${\mathcal D}_d(A_1)\subset{\mathcal D}_d$ the subspace formed by uninodal del Pezzo surfaces.

\begin{proposition}\label{prop1} All but finite number of surfaces in a generic
one-parameter Kodaira-Spencer family
of rational surfaces with $1\le K_\Sigma^2\le 4$ are unnodal
{\rm (}{\it i.e.,} del Pezzo{\rm )},
while the
exceptional members of the family are uninodal del Pezzo.
\end{proposition}

{\bf Proof.}
Let us denote by
${\mathcal T}_{X\Vert D}$ the subsheaf of the sheaf
${\mathcal T}_X$generated by vectors fields tangent to $D$, and by ${\mathcal N}'_{D/X}$
their quotient, so that we obtain the following short exact sequence of sheafs:
$$
0\to{\mathcal T}_{X\Vert D}\to{\mathcal T}_X\to {\mathcal N}'_{D/ X}\to 0.
$$
According to the well known theory of deformations of pairs (see \cite[Section 3.4.4]{Sernesi}),
and due to the long exact cohomology sequence associated to the above short sequence, it is sufficient to show that $h^1(
{\mathcal N}'_{D/X})\ge 2$ if $D$ is either a rational irreducible curve with
$D^2\le - 3$ or $D=D_1
\cup D_2$
where $D_i^2\le - 2$. In the first case, it follows from Serre-Riemann-Roch duality. In the second case,
from the exactness of
the fragment
$H^0({\mathcal N}_{D_2/X})\to H^1({\mathcal N}_{D_1/X})\to H^1({\mathcal N}'_{D/X})\to H^1(
{\mathcal N}_{D_2/X})$
of the long cohomology sequence
associated with the exact sequence of sheaves $0\to {\mathcal N}_{D_1/X}\to {\mathcal N}'_{D/X}\to
{\mathcal N}_{D_2/X}\to 0$.\footnote{In both cases, we use
the equality $H^2({\mathcal T}_{X\Vert D})=0$, which 
can be deduced,
for example,
from Serre duality,
$H^2({\mathcal T}_{X\Vert D})=
(H^0(\Omega^1_X(\log D)\otimes K))^*$,
and Bogomolov-Sommese vanishing $H^0(\Omega^1_X(\log D)\otimes K)=0$;
the latter holds in our case since $X$ is a rational
surface
with $K^2\ge 1$, and thus its anticanonical
Iitaka-Kodaira dimension is equal to $2$.}
\proofend

By a {\it real algebraic surface} we understand a pair
$(Y,c)$, where $Y$ is a complex algebraic surface
and $c:Y\to Y$
is an antiholomorphic involution.
The
classification of minimal real rational surfaces and the classification of real del Pezzo surfaces are
well known:
they are summarized in the two propositions below, respectively
(see, e.g., \cite[Theorems 6.11.11 and 17.3]{DIK}).

\begin{proposition}\label{Co}
Each minimal real rational surface~$Y$ is one of the following:
\begin{enumerate}
\item[(1)]
$\PP^2$ with its standard real structure {\rm(}$d=9${\rm)},
the real part $\R Y$ of $Y$ is homeomorphic to $\R\PP^2$;
\item[(2)]
$\PP^1 \times
\PP^1$ with one of its four nonequivalent real structures
{\rm(}$d=8${\rm ):} $\R Y=(S^1)^2$,
$\R Y=S^2$, and two structures with
$\R Y=\emptyset$;
\item[(3)]
rational geometrically ruled surfaces ${\mathcal F}_a$,
$a\ge2$, with $\R Y=\#_2 \R\PP^2$
and the standard real structure, if
$a$ is odd, and with $\R Y=(S^1)^2$ or~$\emptyset$ and one of the two
respective nonequivalent structures, if $a$ is even {\rm (}$d=8${\rm )};
\item[(4)]
real conic bundles over $
\PP^1$ with $2m\ge4$ reducible fibers, which are
all real and consist of pairs of complex conjugate exceptional curves
{\rm (}$d=8-2m${\rm )}, \quad $\R Y=mS^2$;
\item[(5)]
del
Pezzo surfaces of degree $d=1$ or~$2$:  $\R Y=
\R\PP^2
\sqcup 4S^2$,
if $d=1$, and $\R Y=3S^2$ or $4S^2$, if
$d=2$.
\end{enumerate}
\end{proposition}

\begin{proposition}\label{DDP}
With one exception, a real
del
Pezzo surface~$(Y,c)$ of degree
$d\ge1$ is determined up to deformation by the topology of~$\R Y$.  In
the exceptional case $d=8$ and
$\R Y=\varnothing$, there are two
deformation classes, distinguished by whether $Y/c$ is~$Spin$ or not.

The topological
types of $\R Y$ are the following extremal types and their
derivatives, which are obtained from the extremal ones by sequences
of topological Morse simplifications of $\R Y${\rm:}
\smallskip
\bgroup
\let\comma,
\catcode`\,\active
\def,{\comma\ }
\def\or{\hbox{or\/}\ }
\ialign{\indent\indent$#\:$\ \,\hss&$#$\hss\cr
d=9& \R Y=\R\PP^2;\cr
d=8& \R Y=\#_2 \R\PP^2 \ \or (S^1)^2;\cr
d=7& \R Y=\#_3 \R\PP^2;\cr
d=6& \R Y=\#_4 \R\PP^2 \ \or (S^1)^2;\cr
d=5& \R Y=\#_5 \R\PP^2;\cr
d=4& \R Y=\#_6 \R\PP^2, (S^1)^2,\or 2S^2;\cr
d=3& \R Y=\#_7 \R\PP^2 \ \or \R\PP^2\sqcup S^2;\cr
d=2& \R Y=\#_8 \R\PP^2, 2\R\PP^2, \#_2 \R\PP^2 \sqcup S^2,(S^1)^2,\or 4S^2;\cr
d=1& \R Y=\#_9 \R\PP^2, \#_2 \R\PP^2 \sqcup \R\PP^2, \#_3 \R\PP^2 \sqcup S^2,\or \R\PP^2\sqcup 4S^2.\cr}
\egroup
\end{proposition}

\subsection{Main results}\label{sec1.2}

Let us
consider a real del Pezzo surface $(\Sig,c)$, and
assume that its real point set $\R\Sig=\Fix(c)$
is nonempty. Pick a real divisor class
$D\in\Pic(\Sig)$, satisfying $-DK_\Sig>0$ and $D^2\ge-1$, and put $r=-DK_\Sig-1$.
Fix an integer $m$ such that $0\le 2m\le r$ and
introduce a real structure $c_{r,m}$ on $\Sig^r$ that maps $(w_1,
\ldots , w_r)\in\Sig^r$ to $(w'_1, \ldots , w'_r)\in\Sig^r$ with
$w'_i = c(w_i)$ if $i>2m$, and $(w'_{2j-1}, w'_{2j}) = (c(w_{2j}),
c(w_{2j-1}))$ if $j\le m$. With respect to this real structure a
point $\bw=(w_1,...,w_r)$ is real, {\it i.e.}, $c_{r,m}$-invariant,
if and only if $w_i$ belongs to the real part $\R\Sig$ of $\Sig$ for
$i>2m$ and $w_{2j-1}, w_{2j}$ are conjugate to each other
for $j \le m$. In what follows we work with an open dense subset ${\mathcal
P}_{r,m}(\Sig)$ of $\R\Sig^r=\Fix\, c_{r,m}$ consisting of
$c_{r,m}$-invariant $r$-tuples $\bw=(w_1,...,w_r)$ with pairwise
distinct $w_i\in\Sig$.

Observe that, if
a real irreducible rational curve $C\in|D|$ can be traced through all the
points $w_i$ of $\bw$ and $2m<r=-CK_\Sig-1$, the real
points of $\bw$ must lie
on the unique
one-dimensional
connected component of the real part of $C$,
hence must belong to the same connected component of $\R\Sig$. In
the case $2m=r$,
each real rational curve $C\in|D|$ passing through a collection of
$m$ pairs of complex conjugate points of $\Sig$ has an odd
intersection with the real divisor $K_\Sig$, hence $C$ has a
homologically non-trivial real part in $\R\Sig$.

Thus, we fix a connected component
$F$ of $\R\Sig$ and put
$$
{\mathcal P}_{r,m}(\Sig,F)=\{\bw=(w_1,...,w_r)\in{\mathcal
P}_{r,m}(\Sig)\ :\
w_i\in F\ \text{for}\ i>2m\, \}\ .
$$
Denote by ${\mathcal M}_{0,r}(\Sig,D)$ the set of isomorphism classes of
pairs $(\nu:\PP^1\to\Sig,\bp)$, where $\nu: \PP^1 \to \Sig$ is a holomorphic map
such that $\nu_*(\PP^1)\in|D|$, and $\bp$ is a sequence of $r$ pairwise distinct points in $\PP^1$.
Put
\begin{eqnarray}{\mathcal R}(\Sig,D,F,\bw)&=&\{[\nu:\PP^1\to\Sig,\bp]\in{\mathcal M}_{0,r}(\Sig,D)\ :
\ \nonumber\\
& &\quad \nu\circ\conj=c\circ\nu,\ \nu(\R\PP^1)\subset F,\ \nu(\bp)=\bw\}\ ,\nonumber\end{eqnarray}
where $\conj: \PP^1 \to \PP^1$ is the complex conjugation.

If either the degree of $\Sig$ is greater than $1$, or $D \ne -K_\Sig$, then
for any generic $r$-tuple $\bw\in{\mathcal P}_{r,m}(\Sig,F)$,
the set ${\mathcal R}(\Sig,D,F,\bw)$ is finite and presented by immersions
(see Lemma \ref{ln12}).
In such a case, pick a conjugation-invariant
class \mbox{$\varphi\in H_2(\Sig\setminus F; \Z/2)$} and
put
\begin{equation}W_m(\Sig,D,F,\varphi,
\bw)=\sum_{[\nu,\bp]\in{\mathcal R}(\Sig,D,F,\bw)}
(-1)^{C_+\circ C_- +
C_+
\circ\varphi }\
,\label{enn1}\end{equation} where $C_{\pm}=\nu(\PP^1_\pm)$ with $\PP^1_+,\PP^1_-$ being the two connected components
of $\PP^1\setminus\R\PP^1$.

If the degree of $\Sig$ is equal to $1$ and $D = -K_\Sig$, then
for any generic $r$-tuple $\bw\in{\mathcal P}_{r,m}(\Sig,F)$
and any conjugation-invariant
class \mbox{$\varphi\in H_2(\Sig\setminus F; \Z/2)$}
we define the number $W_m(\Sig, D, F, \varphi, \bw)$
by the formula (\ref{enn1}) retaining in it only the classes $[\nu, \bp]$
presented by immersions.

If $\varphi=0$, we get the original definition of
Welschinger \cite{W,W1}.

\begin{proposition}\label{p1}
The number $W_m(\Sig,D,F,\varphi,\bw)$
does not depend on the choice of a generic element
$\bw
\in {\mathcal
P}_{r,m}(\Sig,F)$.
\end{proposition}

Proposition \ref{p1} is in fact a special case of more general deformation invariance statements.
Consider a smooth real surface $X_0$ with $\R X_0 \ne \varnothing$, a real divisor class $D_0 \in \Pic(X_0)$, a connected component
$F_0$ of $\R X_0$, a conjugation-invariant
class \mbox{$\varphi_0 \in H_2(X_0 \setminus F_0; \Z/2)$}, and a conjugation invariant
collection $\bw_0$ of points in $X_0$.
By an {\it elementary deformation} of the tuple $(X_0, D_0, F_0, \varphi_0)$
(respectively, $(X_0, D_0, F_0, \varphi_0, \bw_0)$)
we mean a one-parameter smooth family of smooth surfaces $X_t$, $t \in [-1, 1]$, extended to a continuous family
of tuples $(X_t, D_t, F_t, \varphi_t)$ (respectively, $(X_t, D_t, F_t, \varphi_t, \bw_t)$).
Two
tuples $T=(X, D, F,\varphi)$ and $T'=(X',D',F',\varphi')$
are called
\emph{deformation equivalent} if they can be connected by a chain
$T=T^{(0)}$, \dots, $T^{(k)}=T'$ so that any two neighboring tuples in the chain are isomorphic
to fibers of an elementary deformation.

\begin{proposition}\label{def-invariance}
Let $(\Sig_t, D_t, F_t, \varphi_t, \bw_t)$, $t\in[-1,1]$, be
an elementary deformation of tuples such that
all surfaces $\Sig_t$,
$t \ne 0$, belong to ${\mathcal D}_d^{\DP}$ for some $1\le d\le 9$,
and the collections $\bw_{\pm 1}$ belong to ${\mathcal P}_{r,m}(\Sig_{\pm 1}, F_{\pm 1})$
and are generic.
Then,
\begin{equation}W_m(\Sig_{-1},D_{-1},F_{-1},\varphi_{-1}, \bw_{-1})=W_m(\Sig_1,D_1,F_1,\varphi_1, \bw_1)\
.\label{enn100}\end{equation}
\end{proposition}

We skip $\bw$ in the notation of the
numbers $W_m(\Sig, D, F, \varphi, \bw)$ and call them
{\it Welschinger invariants}.

Proposition \ref{def-invariance} plays a central role
in the proof of the following statement.

\begin{theorem}\label{t7} If tuples $(\Sig,D,F,\varphi)$ and $(\Sig',D',F',\varphi')$ are deformation equivalent,
then $W_m(\Sig,D,F,\varphi)=W_m(\Sig',D',F',\varphi')$.
\end{theorem}

Proofs of Propositions \ref{p1} and \ref{def-invariance},
as well as the proof of Theorem \ref{t7}, are found in Section \ref{proof}.

\section{Families of rational curves on rational
surfaces}\label{sec-main}

\subsection{General setting}

Let
$\Sig$ be a smooth rational surface,
and $D\in\Pic(\Sig)$
a divisor class. Denote by $\overline{\mathcal M}_{0,n}(\Sig,D)$
the space of the isomorphism classes of
pairs $(\nu:\hat C\to\Sig,\bp)$, where $\hat C$ is either $\PP^1$
or a connected reducible nodal curve of arithmetic genus zero,
$\nu_*\hat C\in|D|$,
$\bp=(p_1,...,p_n)$ is a sequence of distinct smooth points of $\hat C$, and each component of $\hat C$ contracted
by $\nu$
contains at
least three special points. This
moduli space is a projective scheme (see \cite{FP}),
and
there are natural morphisms
$$\Phi_{\Sig,D}:\overline{\mathcal M}_{0,n}(\Sig,D)\to|D|,\quad [\nu:\hat
C\to\Sig,\bp]\mapsto\nu_*\hat C\ ,$$ $$\Ev:\overline{\mathcal
M}_{0,n}(\Sig,D)\to\Sig^n,\quad [\nu:\hat C\to\Sig,\bp]\mapsto\nu(\bp)\
.$$ For any subscheme ${\mathcal
V}\subset\overline{\mathcal M}_{0,n}(\Sig,D)$, define the
{\it intersection dimension} $\idim\mathcal V$ of $\mathcal V$ as follows:
$$\idim{\mathcal V}=\dim(\Phi_{\Sig,D}\times\Ev)({\mathcal V})\ ,$$
where the latter value is the maximum over the dimensions of
all irreducible components.

Put
\begin{eqnarray}
{\mathcal M}_{0,n}^{br}(\Sig,D)&=&\{[\nu:\PP^1\to\Sig,\bp]\ \in {\mathcal M}_{0,n}(\Sig,D)\ :\
\nu\ \text{is birational onto}\ \nu(\PP^1)\},\nonumber\\ {\mathcal M}_{0,n}^{im}(\Sig,D)&=&\{[\nu:\PP^1\to\Sig,\bp]
\in{\mathcal M}_{0,n}(\Sig,D)\ :\ \nu\ \text{is an immersion}\}\ .\nonumber
\end{eqnarray} Denote by $\overline{{\mathcal M}^{br}_{0,n}}(\Sig,D)$ the
closure of ${\mathcal M}^{br}_{0,n}(\Sig,D)$ in
$\overline{\mathcal M}_{0,n}(\Sig,D)$, and introduce also the space
$${\mathcal M}'_{0,n}(\Sig,D)=\{[\nu:\hat C\to\Sig,\bp]
\in
\overline{{\mathcal M}^{br}_{0,n}}(\Sig,D)\ :\ \hat C\simeq\PP^1\}\ .$$

The following 
statement will be used below.

\begin{lemma}\label{ln15}
For any element $[\nu:\PP^1\to\Sig,\bp]\in{\mathcal M}_{0,n}^{br}(\Sig,D)$ such that \mbox{$\nu(\bp)\cap\Sing(\nu(\PP^1))=
\emptyset$}, the map $\Phi_{\Sig,D}\times\Ev$ is
injective
in a neighborhood of
that element,
and, for the germ at $[\nu:\PP^1\to\Sig,\bp]$ of any irreducible subscheme ${\mathcal V}\subset
{\mathcal M}_{0,n}^{br}(\Sig,D)$, we have
$$\dim {\mathcal V}=\idim {\mathcal V}\ .$$
\end{lemma}

{\bf Proof.} The  inequality $\idim {\mathcal V} \le \dim {\mathcal V}$ is immediate from the definition.
The opposite inequality and the injecitvity of  $\Phi_{\Sig,D}\times\Ev$
follow from the observation that, for
an irreducible rational curve $C\in|D|$ and a tuple
$\bz\subset C\setminus\Sing(C)$ of $n$ distinct points, the normalization map
$\nu:\PP^1\to C$ and the lift $\bp=\nu^{-1}(\bz)$ represent the unique preimage of
$(C,\bz)\in|D|\times\Sigma^n$ in ${\mathcal M}_{0,n}^{br}(\Sig,D)$.
\proofend

\subsection{Curves on
del Pezzo and uninodal del Pezzo surfaces}

We
establish here certain properties of the spaces
${\mathcal M}^{im}_{0,0}(\Sig,D)$, ${\mathcal M}^{br}_{0,0}(\Sig,D)$, and $\overline{{\mathcal M}^{br}_{0,0}}(\Sig,D)$,
notably, compute dimension and describe generic members of these spaces as well as of some divisors therein.
These properties
basically
follow
from \cite[Theorem 4.1 and Lemma 4.10]{GP}. However, the cited paper
considers the plane blown up at generic points, whereas we work with arbitrary del Pezzo or uninodal del Pezzo
surfaces.
For this reason, we supply all claims with complete proofs.

Through all this section we use the notation
$$ r=-DK_\Sig-1.$$

\begin{lemma}\label{ln11}
If $\Sig$ is a smooth rational surface and
$-DK_{\Sig}>0$, then
the space ${\mathcal M}^{im}_{0,0}(\Sig,D)$ is either empty, or a smooth
variety of dimension $r$.
\end{lemma}

{\bf Proof}.
Let
$[\nu:\PP^1\to\Sig]\in{\mathcal M}_{0,0}^{im}(\Sig,D)$. The Zariski
tangent space to ${\mathcal M}_{0,0}^{im}(\Sig,D)$ at
$[\nu:\PP^1\to\Sig]$
can be identified with
$H^0(\PP^1,{\mathcal N}^{\nu}_{\PP^1})$, where
${\mathcal N}^{\nu}_{\PP^1}=\nu^*{\mathcal T}\Sig/{\mathcal T}\PP^1$
is the normal bundle. Since
\begin{equation}\deg{\mathcal
N}^{\nu}_{\PP^1}=-DK_\Sig-2\ge-1>(2g-2)\big|_{g=0}=-2\
,\label{en7}\end{equation} we have
\begin{equation}h^1(\PP^1,{\mathcal N}^{\nu}_{\PP^1})=0\ ,\label{e20}\end{equation}
and hence
${\mathcal M}_{0,0}^{im}(\Sig,D)$ is smooth at $[\nu:\PP^1\to\Sig]$
and is
of
dimension
\begin{equation}h^0(\PP^1,{\mathcal N}^{\nu}_{\PP^1})=\deg{\mathcal
N}^{\nu}_{\PP^1}-g+1=-DK_\Sig-1=r\ .\label{en8}\end{equation}
\proofend

\begin{lemma}\label{ln12}
(1) Let $\Sig\in{\mathcal D}_d^{\DP}$ and $-DK_{\Sig}>0$. Then, the following holds:
\begin{enumerate}\item[(i)]
The space ${\mathcal M}_{0,0}^{br}(\Sig,D)$ is either empty or a
variety of dimension $r$, and
$\idim({\mathcal M}_{0,0}(\Sig,D)\setminus{\mathcal M}_{0,0}^{br}(\Sig,D))<r$.
\item[(ii)] If
either $d > 1$ or $D\ne-K_\Sig$, then ${\mathcal
M}_{0,0}^{im}(\Sig,D)\subset{\mathcal M}_{0,0}^{br}(\Sig,D)$ is
an
open dense subset. \item[(iii)]
There
exists
an open dense
subset $U_1
\subset{\mathcal D}^{\DP}_1$ such that, if
$\Sig \in U_1$,
then ${\mathcal M}_{0,0}(\Sig,-K_\Sig)$ consists of
$12$ elements, each corresponding to a rational nodal curve.\end{enumerate}

(2) Let $d\le 4$. There exists an open dense subset $U_d(A_1)\subset{\mathcal D}_d(A_1)$
such that
if $\Sig\in U_d(A_1)$ and $-DK_{\Sig}>0$,
then
\begin{enumerate}
\item[(i)] $\idim{\mathcal M}_{0,0}(\Sig,D)\le r$;
\item[(ii)] a generic
element $[\nu:\PP^1\to\Sig]$ of any irreducible component ${\mathcal V}$ of
${\mathcal M}_{0,0}(\Sig,D)$ such that $\idim{\mathcal V}=r$, is an immersion, and
the divisor
$\nu^*(E_\Sig)$
consists of $DE_\Sig$ distinct points.
\end{enumerate}
\end{lemma}

{\bf Proof}.
Let $\Sig\in{\mathcal D}_d^{\DP}\cup{\mathcal D}_d(A_1)$. All the statements in the case of an effective
$-K_\Sig-D$ immediately follow from elementary properties of plane lines, conics, and cubics.
Thus, we suppose that $-K_\Sig-D$ is not effective.

Let ${\mathcal V}_1$ be an irreducible component of
${\mathcal M}_{0,0}^{br}(\Sig,D)$ and
$[\nu:\PP^1\to\Sig]$ its generic element. Then by \cite[Theorem
II.1.2]{K}
\begin{equation}\dim\Hom(\PP^1,\Sig)_\nu\ge-DK_\Sig+2\chi({\mathcal
O}_{\PP^1})=-DK_\Sig+2\ .\label{en19}\end{equation}
Reducing  by the
automorphisms of $\PP^1$, we get \begin{equation}\dim{\mathcal
V}_1\ge-DK_\Sig+2-3=r\ .\label{en20}\end{equation}
Hence, in
view of Lemma \ref{ln11}, to prove that $\dim{\mathcal
M}_{0,0}^{br}(\Sig,D)=r$ and ${\mathcal M}_{0,0}^{im}(\Sig,D)$ is dense in
${\mathcal M}_{0,0}^{br}(\Sig,D)$, it is enough to show that $\dim ({\mathcal M}_{0,0}^{br}(\Sig,D)\setminus{\mathcal
M}_{0,0}^{im}(\Sig,D))<r$.

Notice, first, that, in the case $r=0$, the
curves $C\in\Phi_{\Sig,D}({\mathcal M}_{0,0}^{br}(\Sig,D))$ are
nonsingular
due to the bound
\begin{equation}-DK_\Sig\ge (C\cdot C')(z)\ge s\ ,\label{en6}\end{equation}
coming from the intersection of $C$ with
a curve $C'\in|-K_\Sig|$ passing through a point $z\in C$, where $C$
has multiplicity $s$. Thus, we suppose that $r>0$.
Let ${\mathcal V}_2$ be an irreducible component
of ${\mathcal M}_{0,0}^{br}(\Sig,D)\setminus {\mathcal
M}_{0,0}^{im}(\Sig,D)$, $[\nu:\PP^1\to\Sig]\in{\mathcal V}_2$ a generic
element, and let $\nu$ have $s\ge1$ critical points
of multiplicities $m_1\ge...\ge m_s\ge2$. In particular, bound (\ref{en6}) gives
\begin{equation}-DK_\Sig\ge m_1\ .\label{en17}\end{equation}
Then ({\it cf.} \cite[First formula in the proof of Corollary
2.4]{CH}),
$$\dim{\mathcal
V}_2\le h^0(\PP^1,{\mathcal
N}^{\nu}_{{\mathcal P}^1}/\Tors({\mathcal N}^{\nu}_{\PP^1}))\ ,$$ where the
normal sheaf ${\mathcal N}^{\nu}_{\PP^1}$ on $\PP^1$ is defined as
the cokernel of the map $d\nu:{\mathcal T}\PP^1\to\nu^*{\mathcal
T}\Sig$, and $\Tors(*)$ is the torsion sheaf.
It follows from \cite[Lemma 2.6]{CH} ({\it cf.} also the computation in \cite[Page 363]{CH})
that $\deg\Tors({\mathcal
N}^{\nu}_{\PP^1})=\sum_i(m_i-1)$, and hence
\begin{equation}\deg{\mathcal N}^{\nu}_{\PP^1}/\Tors({\mathcal
N}^{\nu}_{\PP^1})=-DK_\Sig-2-\sum_{i=1}^s(m_i-1)\label{en24}\end{equation}
which yields
\begin{eqnarray}\dim{\mathcal V}_2
&\le& h^0(\PP^1,{\mathcal
N}^{\nu}_{\PP^1}/\Tors({\mathcal
N}^{\nu}_{\PP^1}))\nonumber\\
&=&\max\{\deg{\mathcal N}^{\nu}_{\PP^1}/\Tors({\mathcal
N}^{\nu}_{\PP^1})+1,\
0\}\overset{\text{(\ref{en17})}}{\le}r-(m_1-1)<r,\label{en5}\end{eqnarray}

Let us show that $\idim{\mathcal V}<r$ for any irreducible component ${\mathcal V}$
of ${\mathcal M}_{0,0}(\Sig,D)\setminus{\mathcal M}_{0,0}^{br}(\Sig,D)$.
Indeed, if a generic element $[\nu:\PP^1\to\Sig]\in{\mathcal V}$ satisfies $\nu_*(\PP^1)=sC$ for some $s\ge2$, then $$\idim{\mathcal
V}\le-\frac{1}{s}DK_\Sig-1<-DK_\Sig-1=r\ .$$

To complete the proof of (2ii), let us assume that $\dim{\mathcal V}=r$ and the divisor $\nu^*(E_\Sig)$
contains a multiple point $sz$, $s\ge2$.
In view of $DE_\Sig \ge s$ and $(-K_\Sig - E_\Sig)D\ge0$
(remind that $D$ is irreducible and $-K-D$ is not effective), we have
$-DK_\Sig\ge s$. Furthermore,  $T_{[\nu]}{\mathcal V}$ can be identified with
a subspace of $H^0(\PP^1,{\mathcal N}_{\PP^1}^{\nu}(-(s-1)z))$ ({\it cf.} \cite[Remark in page 364]{CH}). Since
$$\deg {\mathcal N}_{\PP^1}^{\nu}(-(s-1)z))=-DK_\Sig-1-s\ge-1>-2\ ,$$
we have
$$H^1(\PP^1,{\mathcal N}_{\PP^1}^{\nu}(-(s-1)z))=0\ ,$$ and hence $$
\dim{\mathcal V}\le h^0(\PP^1,{\mathcal N}_{\PP^1}^{\nu}(-(s-1)z))=r-(s-1)<r$$ contrary to the assumption
$\dim{\mathcal V}=r$. \proofend

\begin{lemma}\label{ln17}
There exists an open dense subset $U_2\subset{\mathcal D}_1^{\DP}$
such that, for each $\Sig\in U_2$, the set
of effective divisor classes $D\in\Pic(\Sig)$ satisfying $-DK_\Sig=1$ is finite, the set of rational curves
in the corresponding linear systems $|D|$ is finite, and any two such rational curves $C_1,C_2$ either coincide, or are disjoint, or
intersect in $C_1C_2$ distinct points.
\end{lemma}

{\bf Proof}.
For any $\Sig\in{\mathcal D}^{\DP}_1$, we have $\dim|-2K_\Sig|=3$. Hence, the condition
$-DK_\Sig=1$ yields that $-2K_\Sig-D$ is effective, which in turn implies the finiteness of
the set of effective divisors such that $-DK_\Sig=1$. The finiteness of the set of rational curves in these linear systems
$|D|$ follows from Lemma \ref{ln12}(i). At last, for a generic $\Sig\in{\mathcal D}_1^{\DP}$, these curves are either
singular elements in the elliptic pencil $|-K_\Sig|$ or the $(-1)$-curves,
and as it follows easily, for example, from considering $\Sig$ as a projective plane blown up
at $8$ generic points, any two
of these curves
intersect transversally and in distinct smooth points.
\proofend

\begin{lemma}\label{ln16}
Let $U_1$, $U_2$ be the subsets of  ${\mathcal D}_1^{\DP}$ introduced in Lemmas \ref{ln12} and \ref{ln17}, respectively.
For each $\Sig\in U_1\cap U_2$,
each $D \in \Pic(\Sig)$ with
$-DK_\Sig>0$ and $D^2\ge -1$,
and for each irreducible component ${\mathcal V}$ of
$\overline{{\mathcal M}^{br}_{0,0}}(\Sig,D)\setminus{\mathcal M}_{0,0}^{br}(\Sig,D)$
with $\idim{\mathcal V}=r-1$,
one has:
\begin{enumerate}
\item[(i)] A generic element $[\nu:\hat C\to\Sig]\in{\mathcal V}$ is as follows
\begin{itemize}\item
$\hat C=\hat C_1\cup\hat C_2$ with $\hat C_i\simeq\PP^1$,
$[\nu
|_{\hat C_i}:\hat C_i\to\Sig]\in{\mathcal M}_{0,0}^{im}(\Sig,D_i)$,
where $D_1D_2>0$ and $-D_iK_\Sig>0$, $D_i^2\ge -1$
for each $i=1,2$;
\item $\nu(\hat C_1) \ne \nu(\hat C_2)$, except for the only case when
$D_1=D_2=-K_\Sig$ and $\nu(\hat C_1) = \nu(\hat C_2)$ is one of the $12$ uninodal curves in $|-K_\Sig|$;
\item
$\nu$ is an immersion {\rm (}{\it i.e.}, a local isomorphism onto the image{\rm )}.
\end{itemize}
Moreover, each element $[\nu:\hat C\to\Sig]\in\overline{\mathcal M}_{0,0}(\Sig,D)$
as above does belong to $\overline{{\mathcal M}^{br}_{0,0}}(\Sig,D)$.
\item[(ii)] The germ of $\overline{{\mathcal
M}^{br}_{0,0}}(\Sig,D)$ at a generic element of ${\mathcal V}$ is
smooth.\end{enumerate}
\end{lemma}

{\bf Proof}.
Show, first, that $\idim({\mathcal M}'_{0,0}(\Sig,D)\setminus{\mathcal
M}_{0,0}^{br}(\Sig,D))\le r-2$. Assume on the contrary that there exists
a component ${\mathcal V}$ of ${\mathcal
M}'_{0,0}(\Sig,D)\setminus{\mathcal M}_{0,0}^{br}(\Sig,D)$ with
$\idim{\mathcal V}=r-1$ ($\idim{\mathcal V}$ cannot be bigger by Lemma
\ref{ln12}(i)).
Then its generic element $[\nu:\PP^1\to\Sig]$ is such that $\nu_*(\PP^1)=sC$ with $C$ an irreducible rational curve, $s\ge2$.
Thus, $$r-1=-sCK_\Sig-2\le-CK_\Sig-1=\dim{\mathcal M}_{0,0}^{br}(\Sig,C)\ ,$$
which yields $s=2$ and $-CK_\Sig=1$. By adjunction formula, either
$C^2=-1$, or $C^2\ge1$. The former case is excluded by the assumption $D^2\ge-1$.
In the case $C^2\ge1$, since $K_\Sig^2=1$ and $-CK_\Sig=1$, the only possibility is $C\in|-K_\Sig|$.
However, in such a case
the map $\nu$ cannot be deformed into an element of ${\mathcal M}^{br}_{0,0}(\Sig,-2K_\Sig)$, since
$C$ has a node, and hence the
deformed map would birationally send $\PP^1$ onto a curve with
$\delta$-invariant $\ge 4$,
which is bigger
than its arithmetic genus, $((-2K_\Sig)^2+(-2K_\Sig)K_\Sig)/2+1=2$.

Let $[\nu:\hat C\to\Sig]$ be a generic element of an irreducible component
${\mathcal V}$ of $\overline{{\mathcal M}^{br}_{0,0}}(\Sig,D)\setminus{\mathcal
M}'_{0,0}(\Sig,D)$ with $\idim{\mathcal V}=r-1$. Then $\hat C$ has at least $2$ components.
On the other side, if $\hat C$ had $\ge3$ components, Lemma \ref{ln12}(1) would yield $\idim{\mathcal V}
\le-DK_\Sig-3<r-1$. Hence
$\hat C=\hat C_1\cup\hat C_2$, $\hat C_1\simeq\hat C_2\simeq\PP^1$,
and, according to Lemma \ref{ln11} and Lemma \ref{ln12}(1),
for each $i=1,2$ we have: $\nu_i=\nu\big|_{\hat C_i}$ is an immersion,
$\dim{\mathcal M}_{0,0}(\Sig,D_i)_{[\nu_i]}=-D_iK_\Sig-1$, and $-D_iK _\Sig>0$, $D_i^2\ge-1$.

If $-DK_\Sig=2$ and $\nu(\hat C_1)\ne\nu(\hat C_2)$, then the intersection points of these curves
are nodes, which follows from the definition of the set $U_2$ (see Lemma \ref{ln17}), and hence
$\nu$ is an immersion at the node $\hat z$ of $\hat C$.

If $-DK_\Sig=2$ and $\nu(\hat C_1)=\nu(\hat C_2)$, then $D_1=D_2$ and $D_1^2=D_2^2\ge1$
in view of the adjunction formula and the condition $D^2\ge-1$. It is easy to see that this is only possible,
when $D_1=D_2=-K_\Sig$. In particular, by the definition of the set $U_1$ (see Lemma \ref{ln12}(iii)), the curve
$C=\nu(\hat C_1)=\nu(\hat C_2)\in|-K_\Sig|$ has one node $z$.
We then see that, $\nu$ takes the germ $(\hat C,\hat z)$ isomorphically onto the germ
$(C,z)$, since, otherwise we would get a deformed map $\nu$ with
the image whose
$\delta$-invariant $\ge 4$,
which is bigger
than its arithmetic genus, $((-2K_\Sig)^2+(-2K_\Sig)K_\Sig)/2+1=2$.

Suppose, now, that $-DK_\Sig>2$, thus, $-D_1K_\Sig>1$. Then $$\dim{\mathcal M}_{0,0}(\Sig,D_1)_{[\nu_1]}>0\ ,$$
and hence $C_1\ne C_2$. To prove that $\nu$ is an immersion at the node $\hat z\in\hat C$, we will show that
any two local branches of $\nu_1$ and $\nu_2$ either are disjoint, or intersect transversally. Indeed, assume
on the contrary that
there exist $z_i\in\hat C_i$, $i=1,2$, such that $\nu_1(z_1)=\nu_2(z_2)=z\in\Sig$,
and $\nu_1(\hat C_1,z_1)$ intersects $\nu_2(\hat C_2,z_2)$ at $z$
with multiplicity $\ge2$.
Then
\begin{equation}\dim{\mathcal M}_{0,0}(\Sig,D_1)_{[\nu_1]}\le h^0(\hat C_1,{\mathcal N}_{\hat C_1}^{\nu_1}(-z_1))\ .\label{e13}\end{equation}
Since
$$\deg{\mathcal N}_{\hat C_1}^{\nu_1}(-z_1)=-D_1K_\Sig-2-1=-D_1K_\Sig-3
>-2\ ,$$
we get
$h^1(\hat C_1,{\mathcal N}_{\hat C_1}^{\nu_1}(-z_1))=0$.
Therefore,
$$\deg{\mathcal N}_{\hat C_1}^{\nu_1}(-z_1)\le \deg{\mathcal N}_{\hat C_1}^{\nu_1}(-z_1)+1=-D_1K_\Sig-2<-D_1K_\Sig-1
=\dim{\mathcal M}_{0,0}(\Sig,D_1)_{[\nu_1]}\ ,$$
which contradicts (\ref{e13}).

The smoothness of $\overline{{\mathcal M}^{br}_{0,0}}(\Sig,D)$ at $[\nu:\hat
C\to\Sig]$, where $\nu_*\hat C$ is a reduced nodal curve, follows from \cite[Lemma 2.9]{T}, where the requirements
are $D_iK_\Sig<0$, $i=1,2$. We
will
show that the same requirements suffice under assumption that $\nu$ is an immersion.
Let us show that
\begin{equation}T_{[\nu]}\overline{{\mathcal M}^{br}_{0,0}}(\Sig,D)\simeq H^0(\hat C,{\mathcal N}_{\hat C}^\nu)\ ,
\label{eK4}\end{equation} where the normal sheaf ${\mathcal N}_{\hat C}^\nu$
comes from the exact sequence
\begin{equation}0\to{\mathcal T}_{\hat C}\to\nu^*{\mathcal T}_\Sig\to{\mathcal N}_{\hat C}^\nu\to0\ ,
\label{eK1}\end{equation}
${\mathcal T}_\Sig$ being the tangent bundle of $\Sig$, and ${\mathcal T}_{\hat C}$ being the
tangent sheaf of $\hat C$ viewed as the push-forward by the normalization
$\pi:\hat C_1\sqcup\hat C_2\to\hat C$ of the subsheaf ${\mathcal T}'_{\hat C_1\sqcup\hat C_2}\subset{\mathcal T}_{\hat C_1\sqcup\hat C_2}$
generated by the sections vanishing at the preimages of the node $z\in\hat C$.

Indeed,
the Zariski tangent space to $\Hom(\hat C,\Sig)$ at $\nu$ is naturally isomorphic to
$H^0(\hat C,\nu^*{\mathcal T}_\Sig)$ (see \cite[Theorem 1.7, Section II.1]{K}). Next, we take the quotient by action of the germ
of $\Aut(\hat C)$ at the identity. This germ is smooth and acts freely on the germ of $\Hom(\hat C,\Sig)$ at $\nu$.
The tangent space
to $\Aut(\hat C)$ at the identity is isomorphic to $H^0(\hat C,{\mathcal T}_{\hat C})$ ({\it cf.} \cite[2.16.4, Section I.2]{K}).
Since
\begin{equation}H^1(\hat C,{\mathcal T}_{\hat C})=
H^1(\hat C_1\sqcup\hat C_2,{\mathcal T}'_{\hat C_1\sqcup\hat C_2})=H^1(\hat C_1,{\mathcal O}_{\hat C_1}(1))\oplus H^1(\hat C_2,{\mathcal O}_{\hat C_2}(1))=0\ ,
\label{eK3}\end{equation}
the associated to (\ref{eK1}) cohomology exact sequence yields
$$T_{[\nu]}\overline{{\mathcal M}^{br}_{0,0}}(\Sig,D)\simeq T_\nu\Hom(\hat C,\Sig)/T_{\Id}\Aut(\hat C)\simeq
H^0(\hat C,\nu^*{\mathcal T}_\Sig)/H^0(\hat C,{\mathcal T}_{\hat C})\simeq
H^0(\hat C,{\mathcal N}_{\hat C}^\nu)\ .$$
We will verify that
\begin{equation}h^0(\hat C,{\mathcal N}_{\hat C}^\nu)=r\ ,
\label{eK2}\end{equation} which in view of
$\dim_{[\nu]}\overline{{\mathcal M}^{br}_{0,0}}(\Sig,D)=r$ (see Lemma \ref{ln12}(i)) will imply the
smoothness of $\overline{{\mathcal M}^{br}_{0,0}}(\Sig,D)$ at $[\nu]$.
There exists a natural morphism of sheaves on $\hat C$:
$$\alpha:\pi_*{\mathcal
N}^{\nu\circ\pi}_{\hat C_1\sqcup
\hat C_2}{\longrightarrow}{\mathcal
N}^{\nu}_{\hat C}\ ,$$ where
$\alpha$ is an isomorphism outside $z$ and acts at $z$ as
follows: since $\nu$ embeds the germ of $\hat C$ at $z$ into $\Sig$,
one can identify the stalk $\left(\pi_*{\mathcal
N}^{\nu\circ\pi}_{\hat C_1\sqcup\hat C_2}\right)_{z}$ with
$\C\{x\}\oplus\C\{y\}$, the stalk $({\mathcal N}_{\hat C}^\nu)_z$
with $\C\{x,y\}/\langle xy\rangle$, and write
$$\alpha_{z}(f(x),g(y))=xf(x)+yg(y)\in\left({\mathcal
N}^{\nu}_{\hat C}\right)_{z}\cong\C\{x,y\}/\langle xy\rangle\ .$$ Hence we obtain an exact sequence of sheaves
\begin{equation}0\to \pi_*{\mathcal
N}^{\nu\circ\pi}_{\hat C_1\sqcup
\hat C_2}\overset{\alpha}{\longrightarrow}{\mathcal
N}^{\nu}_{\hat C}\to{\mathcal O}_{z}\to0\
,\label{en30}\end{equation} whose cohomology sequence vanishes at
$$h^1(z,{\mathcal O}_{z})=0,\quad
h^1(\hat C_1\sqcup\hat C_2,{\mathcal
N}^{\nu\circ\pi}_{\hat C_1\sqcup\hat C_2})=0\ ,$$ (the latter one is equivalent to (\ref{e20}));
hence $h^1(\hat C,{\mathcal N}_{\hat C}^\nu)=0$ and, furthermore,
$$h^0(\hat C,{\mathcal N}_{\hat C}^\nu)=h^0(\hat C_1\sqcup\hat C_2,{\mathcal N}_{\hat C_1\sqcup\hat C_2}^{\nu\circ\pi}))+
h^0(z,{\mathcal O}_z)=
h^0(\hat C_1,{\mathcal N}_{\hat C_1}^{\nu_1})+h^0(\hat C_2,{\mathcal N}_{\hat C_2}^{\nu_2})
+h^0(z,{\mathcal O}_z)$$
$$\overset{\text{{\it cf.}\ (\ref{en8})}}{=}(-D_1K_\Sig-1)+(-D_2K_\Sig-1)+1=r$$
as predicted in (\ref{eK1}).

Finally, let us show that any element $[\nu:\hat C\to\Sig]\in\overline{\mathcal M}_{0,0}(\Sig,D)$,
satisfying conditions of (1i-iii), belongs to $\overline{{\mathcal M}^{br}_{0,0}}(\Sig,D)$,
or, equivalently, admits a deformation into a map $\PP^1\to\Sig$ birational onto its image. Indeed, it follows from
\cite[Theorem 15]{AK} under the condition $h^1(\hat C,\nu^*{\mathcal T}_\Sig)=0$, which one
obtains from the cohomology exact sequence associated with (\ref{eK1}) and
from vanishing relations (\ref{eK3}) and (\ref{eK2}).\proofend

\begin{lemma}\label{p2}
Consider the subsets $U_1$, $U_2$ of  ${\mathcal D}_1^{\DP}$ introduced in Lemmas \ref{ln12} and \ref{ln17}, respectively,
a surface $\Sig\in U_1\cap U_2\subset{\mathcal D}_1^{\DP}$, and
an effective divisor class $D\in\Pic(\Sig)$ such that $-DK_\Sig\ge2$.
Let
$\bw=(w_1,...,w_r)$ be a sequence of $r$ distinct points in $\Sig$, let
$\sig_i$
be
smooth curve germs in $\Sig$ centered at $w_i$, $r'<i\le
r$, for some $r'<r$, $\bw'=(w_i)_{1\le i\le r'}$, and let
\begin{eqnarray}\overline{{\mathcal M}_{0,r}^{br}}(\Sig,D;\bw',\{\sig_i\}_{r'<i\le
r})&=&\{[\nu:\hat C\to\Sig,\bp]\in\overline{{\mathcal
M}_{0,r}^{br}}(\Sig,D)\ : \nonumber\\
& & \quad \nu(p_i)=w_i \,\text{for}\, \ 1\le i\le r',\ \nu(p_i)\in \sig_i, \,\text{for}\, \ r'<i\le
r\}\ .\nonumber\end{eqnarray}

(1) Suppose that $[\nu:\PP^1\to\Sig,\bp]\in\overline{{\mathcal
M}_{0,r}^{br}}(\Sig,D;\bw)\cap{\mathcal
M}_{0,r}^{im}(\Sig,D)$. Then $\Ev$ sends the germ of $\overline{{\mathcal
M}_{0,r}^{br}}(\Sig,D;\bw',\{\sig_i\}_{r'<i\le r})$ at
$[\nu:\PP^1\to\Sig,\bp]$ diffeomorphically onto $\prod_{r'<i\le
r}\sig_i$.

(2) Suppose that $[\nu:\hat C\to\Sig,\bp]\in\overline{{\mathcal
M}_{0,r}^{br}}(\Sig,D;\bw)$ is such that
\begin{itemize}\item $[\nu:\hat C\to\Sig]\in\overline{{\mathcal
M}_{0,0}^{br}}(\Sig,D)$ is as in Lemma \ref{ln16}(i),
\item $r'\ge-D_1K_\Sig-1$, $\#(\bp\cap \hat
C_1)=-D_1K_\Sig-1$, $\#(\bp\cap \hat C_2)=-D_2K_\Sig$, the point sequences
$(w_i)_{1\le i<-D_1K_\Sig}$, $(w_i)_{-D_1K_\Sig\le i\le r}$ are
generic on $C_1=\nu_*\hat C_1$, $C_2=\nu_*\hat
C_2$, respectively, and the germs $\sig_i$, $r'<i\le r$, cross $C_2$
transversally.\end{itemize} Then $\Ev$ sends the germ of
$\overline{{\mathcal
M}_{0,r}^{br}}(\Sig,D;\bw',\{\sig_i\}_{r'<i\le r})$ at
$[\nu:\hat C\to\Sig,\bp]$ diffeomorphically onto $\prod_{r'<i\le
r}\sig_i$.
\end{lemma}

{\bf Proof.} Both statements follow from the fact that $\Ev$ diffeomorphically sends the
germ of $\overline{{\mathcal M}_{0,r}^{br}}(\Sig,D)$ at $[\nu:\hat C\to\Sig,\bp]$ onto the germ
of $\Sig^r$ at $\bw=\nu(\bp)$.

In view of $\dim\overline{{\mathcal M}_{0,r}^{br}}(\Sig,D)=2r$
Lemma \ref{ln12}(i)),
it is
sufficient
to show that
the Zariski tangent space to $\Ev^{-1}(\bw)$ is zero-dimensional.
In view of relation (\ref{eK4}) this is
equivalent to
\begin{equation}h^0(\hat C,{\mathcal N}_{\hat C}^\nu(-\bp))=0\ .\label{eK5}\end{equation}

In the case of $[\nu:\PP^1\to\Sig,\bp]\in\overline{{\mathcal
M}_{0,r}^{br}}(\Sig,D;\bw',\{\sig_i\}_{r'<i\le r})\cap{\mathcal
M}_{0,r}^{im}(\Sig,D)$, we have
$$\deg{\mathcal N}_{\hat C}^\nu(-\bp)=(-DK_\Sig-2)-(-DK_\Sig-1)=-1>-2\ ,$$ and hence (\ref{eK5}) follows
by Riemann-Roch.

In the second case, put $\widetilde\bp=\bp\setminus\{p_r\}$
and twist the
exact sequence (\ref{en30}) to get
$$0\to\pi_*{\mathcal
N}^{\nu\circ\pi}_{\hat C_1\sqcup
\hat C_2}(-\widetilde\bp)\to{\mathcal
N}^{\nu}_{\hat C}(-\widetilde\bp)\to{\mathcal O}_{z}\to0.
$$
Since
$$\deg{\mathcal N}_{\hat C_i}^{\nu_1}(-\widetilde\bp\cap\hat C_i)=(-D_iK_\Sig-2)-(-D_iK_\Sig-1)=-1>-2,
\quad i=1,2\ ,$$
we have
$h^1(\pi_*{\mathcal
N}^{\nu\circ\pi}_{\hat C_1\sqcup
\hat C_2}(-\widetilde\bp))=0$,
and
$h^0(\hat C,\pi_*{\mathcal
N}^{\nu\circ\pi}_{\hat C_1\sqcup
\hat C_2}(-\widetilde\bp))=0$, which
yields that
$H^0(\hat C,{\mathcal
N}^{\nu}_{\hat C}(-\widetilde\bp))$ is isomorphically mapped onto
$H^0(z,{\mathcal O}_z)\simeq\C$. It
implies
that a non-zero global
section of the sheaf ${\mathcal
N}^{\nu}_{\hat C}(-\widetilde\bp)$ does not vanish at $z$, and hence, it does not vanish at $p_r$ chosen on $\hat C_2$ in a generic
way. Thus, (\ref{eK5}) follows. \proofend

\subsection{Deformation of isolated curve singularities}

Let
us recall a few facts on deformations of
curve
singularities
(see, for example, \cite{DH}).
Let $\Sig$ be a smooth algebraic surface,
$z$ an isolated singular point of a curve $C\subset\Sig$, and $B_{C,z}$
the base of a semiuniversal
deformation
of the germ $(C,z)$.
This base
can
be viewed as
a
germ $(\C^N,0)$
and can
be identified with ${\mathcal O}_{C,z}/J_{C,z}$, where
$J_{C,z}\subset{\mathcal O}_{C,z}$
is the Jacobian ideal.

The
equigeneric locus
$B^{\,eg}_{C,z} \subset B_{C,z}$
parametrizes local
deformations of $(C,z)$ with constant
$\delta$-invariant equal to $\del(C,z)$.
This
locus
is irreducible
and has codimension
$\del(C,z)$ in $B_{C,
z}$.
The subset $B^{\,eg,im}_{C,z} \subset B^{\,eg}_{C,z}$
that parametrizes the immersed deformations is open and dense in  $B^{\,eg}_{C, z}$,
and consists only of smooth points of $B^{\,eg}_{C, z}$.
The tangent cone
$T_0B^{\,eg}_{C,z}$ (defined as the limit of the
tangent spaces at points of $B^{\,eg,im}_{C,z}$) can be
identified with $J^{cond}_{C,z}/J_{C,z}$, where \mbox{$J^{cond}_{C, z}\subset{\mathcal O}_{C,z}$}
is the conductor
ideal.
The subset $B^{\,eg,nod}_{C,z} \subset B^{\,eg}_{C,z}$ that
parameterizes the
nodal deformations is also open and dense.
Furthermore, $B^{\,eg}_{C,z}\setminus B^{\,eg,nod}_{C,z}$ is the closure of
three codimension-one strata:
$B^{\,eg}_{C,z}(A_2)$
that parameterizes deformations with one cusp
$A_2$ and $\del(C,z)-1$ nodes,
$B^{\,eg}_{C,z}(A_3)$
that parameterizes
deformations with one tacnode $A_3$ and $\del(C,z)-2$ nodes, and $B^{\,eg}_{C,
z}(D_4)$
that parameterizes deformations with one ordinary triple point
$D_4$ and $\del(C,z)-3$ nodes.

If $C \subset \Sig$ is a curve with isolated singularities,
we consider the joint semiuniversal deformation
for all singular points of $C$. The
base of this deformation,
the equigeneric locus,
and the tangent cone to this locus at the point corresponding to $C$ are as follows:
$$
B_C=\prod_{z\in\Sing(C)}B_{C,z},\quad
B^{\,eg}_C=\prod_{z \in \Sing(C)}B^{\,eg}_{C,z}, \quad
T_0B^{\,eg}_C=\prod_{z \in \Sing(C)}T_0B^{\,eg}_{C,z}\ .$$

\begin{lemma}\label{leg} Let $\nu:\PP^1\to\Sig$ be birational onto its image $C=\nu(\PP^1)$.
Assume that $C\in|D|$, where $D$ is a divisor class such that $r=-DK_\Sig-1>0$.
Let $\bp$ be an $r$-tuple of distinct points of $\PP^1$ such that
$\bw=\nu(\bp)$ is an $r$-tuple of distinct nonsingular points of $C$.
Let $|D|_\bw\subset|D|$
be the linear subsystem of curves passing through
$\bw$, and $\Lam(\bw)\subset B_C$ be the natural image of $|D|_\bw$.
\begin{enumerate}
\item[(1)] One has $\codim_{B_C}\Lam(\bw)=\dim B^{\,eg}_C$, and $\Lam(\bw)$ intersects $T_0B^{\,eg}_C$
transvesally.
\item[(2)] For any $r$-tuple $\widetilde\bw \in \Sig^r$ sufficiently close to $\bw$ and such that
$\Lam(\widetilde\bw)$ intersects $B^{\,eg}_C$
transversally and only at smooth points,
the natural map from the germ ${\mathcal
M}_{0,r}(\Sig,D)_{[\nu,\bp]}$ of ${\mathcal M}_{0,r}(\Sig,D)$ at $[\nu:\PP^1\to\Sig,\bp]$
to $B^{\,eg}_C$ gives rise to a bijection
between the set of elements $[\widetilde\nu:\PP^1\to\Sig,\widetilde\bp]
\in{\mathcal M}_{0,r}(\Sig,D)_{[\nu,\bp]}$ such that $\widetilde\nu(\widetilde\bp) = \widetilde\bw$
on one side
and the set
$\Lam(\widetilde\bw)\cap B^{\,eg}_C$ on the other side.
\end{enumerate}
\end{lemma}

{\bf Proof}.
(1)
The dimension and the transversality statements reduce to the fact that the pull-back of $T_0B^{\,eg}_C$ to $|D|$
intersects
$|D|_\bw$  transversally
and only at one point. In view of
the identification of $T_0B^{\,eg}_C$ with $\prod_{z \in \Sing(C)}J^{cond}_{C,z}/J_{C,z}$
\cite[Theorem 4.15]{DH}, both
required claims read as
\begin{equation}
H^0(C,{\mathcal J}^{cond}_C(-\bw)\otimes{\mathcal O}_\Sig(D))=0\ ,\label{eeg}\end{equation} where
${\mathcal J}^{cond}_C=\Ann(\nu_*{\mathcal O}_{\PP^1}/{\mathcal O}_C)$ is the conductor ideal sheaf,
since ${\mathcal J}^{cond}_C$ can be equivalently regarded as the ideal sheaf of
the zero-dimensional subscheme of $C$
defined at all singular points $z\in\Sing(C)$ by the conductor ideals $J^{cond}_{C,z}=\Ann(\nu_*
\bigoplus_{q\in\nu^{-1}(z)}{\mathcal O}_{\PP^1,q})/{\mathcal O}_{C,z}$.
It is known that
${\mathcal J}^{cond}_C=\nu_*{\mathcal O}_{\PP^1}(-\Delta)$, where $\Delta\subset\PP^1$ is the so-called double-point divisor,
whose degree is $\deg\Delta=2\sum_{z\in\Sing(C)}\del(C,z)$ (see, for example, \cite[Section 2.4]{CH} or
\cite[Section 4.2.4]{DI}). Hence, the relations (\ref{eeg}) can be rewritten as
\begin{equation}
H^0(\PP^1,{\mathcal O}_{\PP^1}(\bd-\Delta-\bp))=0\ ,\label{eeg1}\end{equation} where $\deg\bd=D^2$. Since
$$
\displaylines{
\deg{\mathcal O}_{\PP^1}(\bd-\Delta-\bp))=D^2-2\sum_{z\in\Sing(C)}\del(C,z)-r \cr
= D^2-2\left(\frac{D^2+DK_\Sig}{2}+1\right)-(-DK_\Sig-1)=-1>-2\ ,
}
$$ we obtain $H^1(\PP^1,{\mathcal O}_{\PP^1}(\bd-\Delta-\bp))=0$, and hence by Riemann-Roch
$$\dim H^0(\PP^1,{\mathcal O}_{\PP^1}(\bd-\Delta-\bp))=\deg{\mathcal O}_{\PP^1}(\bd-\Delta-\bp))+1=0\ .$$

(2) The second statement of Lemma immediately follows from the first one due to
the fact that the tangent spaces to the stratum $B^{\,eg}_C$ at its smooth points close to
the origin converge to the same linear space of dimension
$\dim B^{\,eg}_C$ \cite[Theorem 4.15]{DH}.
\proofend

Suppose now that $\Sig$ possesses a real structure, $C$ is a real curve, and $z$ is its real singular point.
Let $b\in B^{\,eg,im}_{C,z}$ be a real point, and let $C_b$
be the corresponding fiber
of the semiuniversal deformation of the germ $(C,z)$.
Define the
Welschinger sign $W_b$ as follows. Let $\pi:\hat C_b\to C_b\hookrightarrow\Sig$ be
the normalization of $C_b$. Here $\hat C_b$ is the union of discs, some of them being real
({\it i.e.}, invariant with respect to the complex conjugation),
the others
forming complex conjugate pairs. Put
$W_b
=(-1)^{C_{b,+}\circ C_{b,-}}$,
where
$C_{b,\pm}=\pi(\hat C_{b,\pm})$ and $\hat C_b\setminus\R\hat C_b=\hat C_{b,+} \sqcup \hat C_{b,-}$ is a
splitting into disjoint complex conjugate halves.

\begin{lemma}\label{new_lemma1}
The Welschinger sign $W_b$ is equal to $(-1)^s$, where $s$ is the number of solitary
nodes in a small real nodal
perturbation
of $C_b$.
\end{lemma}

{\bf Proof}. Straightforward from the definition.
\proofend

\begin{lemma}\label{ln2}
Let $L_t$, $t\in(-\eps,\eps)\subset\R$, be a
smooth one-parameter family of conjugation-invariant affine
subspaces of $B_{C,z}$ of dimension $\del(C,z)$ such that
\begin{itemize}\item $L_0$ passes through the origin and is
transversal to $T_0B^{\,eg}_{C,z}$, \item $L_t\cap B^{\,eg}_{C,z}\subset
B^{\,eg,im}_{C,z}$ for each $t \in (-\varepsilon, \varepsilon) \setminus \{0\}$.
\end{itemize} Then,
\begin{enumerate}
\item[(i)]
the intersection $L_t\cap B^{\,eg}_{C,z}$ is finite for each $t \in (-\varepsilon', \varepsilon') \setminus \{0\}$,
where $\varepsilon' > 0$ is sufficiently small.
\item[(ii)] the function $W(t)=\sum_{b\in L_t\cap\R B^{\,eg}_{C,z}}W_b$ is
constant in $(-\epsilon', \epsilon') \setminus \{0\}$, where $\varepsilon' > 0$ is sufficiently small.
\end{enumerate}
\end{lemma}

{\bf Proof.} The finiteness of the intersection follows from the
transversality of
$L_0$ and $T_0B^{\,eg}_{C,z}$ in $B_{C,z}$.
To prove the second statement, assume, first, that the germ $(C,z)$ represents an ordinary cusp
$A_2$. Then $\R B_{C,z}=(\R^2,0)$ and $\R B^{\,eg}_{C,z}$ is a semicubical
parabola with vertex at the origin.
For the points $b$ belonging to
one of the two connected components of $\R B^{\,eg}_{C,z}\setminus\{0\}$,
the curve $C_b$ has a non-solitary real node;
for the points $b$
from the other component, $C_b$ has a solitary node.
Since, in addition, the line $L_0$ crosses the tangent to
the parabola at the origin transversally
we have $W(t)=0$ for
each $t \in (-\varepsilon', \varepsilon') \setminus \{0\}$ for sufficiently small $\varepsilon' > 0$.

In the general case,
if $\varepsilon' > 0$ is sufficiently small, then for any two points
$t_1 < t_2$ in $(-\epsilon', \epsilon') \setminus \{0\}$ we can connect $L_{t_1}$ with $L_{t_2}$
by a family of $\del(C, z)$-dimensional conjugation-invariant affine subspaces
$L'_t \subset B_{C, z}$, $t\in[t_1,t_2]$, such that
\begin{itemize}
\item the subspaces $L'_t$, $t \in [t_1, t_2]$,
are transversal to $B^{\,eg}_{C,z}$,
\item the intersection number of $L'_t$ and $B^{\,eg}_{C,z}$
is constant in $[t_1, t_2]$,
\item for all but finitely many values of $t$
the intersection $L'_t \cap B^{\,eg}_{C,z}$ is contained in $B^{\,eg,nod}_{C,z}$,
and for the remaining values of $t$, the subspace $L'_t$ intersects $\R B^{\,eg}_{C,z}$ within
$B^{\,eg}_{C,z}(A_2)\cup B^{\,eg}_{C,z}(A_3)\cup B^{\,eg}_{C,z}(D_4)$.
\end{itemize}
The bifurcations through the immersed singularities $A_3$ and $D_4$
do not affect $W(t)$, as well as the cuspidal bifurcation, which we have treated above.
\proofend

\begin{remark}
In fact, Lemma \ref{ln2}
allows one to extend the definition of Welschinger signs
and attribute a {\it Welschinger weight} to any map
$\nu:\PP^1\to\Sigma$ birational onto its image.
\end{remark}

\section{Proof of
Theorem \ref{t7}}\label{proof}
\subsection{Preliminary observations}
We start with two
remarks.

(1) If $Y$ is an irreducible complex variety, equipped with a real structure, and $\R Y$ contains
nonsingular points of $Y$, then $\R Y\cap U\ne\emptyset$ for any Zariski open subset $U\subset Y$.
In particular, a generic element of ${\mathcal P}_{r,m}(\Sig,F)$ is generic in $\Sig^r$.

(2) By blowing up extra real points we can reduce the consideration to the case of del Pezzo surfaces of
degree $1$.

The following statement will be used in the sequel.

\begin{lemma}\label{l45}
Let $t\in(\R,0)\mapsto\Sig_t$ be
a germ of
an elementary deformation $(\Sig_t,D_t,F_t,\varphi_t,\bw_t)$
of a tuple $(\Sig_0,D_0,F_0,\varphi_0,\bw_0)$, where
$\Sig_0$ is a del Pezzo surface of degree $1$,
$D_0\in\Pic(\Sig_0)$ is a real effective divisor
such that $r=-D_0K_{\Sig_0}-1>0$,
and $\bw_0$ belongs to ${\mathcal P}_{r,m}(\Sig_0,F_0)$
and is generic.
Then
$$W_m(\Sig_t,D_t,F_t,\varphi_t,\bw_t)=W_m(\Sig_0,D_0,F_0,\varphi_0,\bw_0)\ .$$
\end{lemma}

{\bf Proof.}
Since $D_0K_{\Sig_0} > 1$ and $\bw_0$ is generic, Lemma \ref{ln12} implies that
all the curves under count are immersed. Thus, each of these curves contributes $1$
to the Gromov-Witten invariant, and the required equality
follows from
Lemma \ref{new_lemma1}.
\proofend

\subsection{Proof of Proposition \ref{p1}}
The only
situation
to consider is the one
where $\Sig \in {\mathcal D}_1^{\DP}$ and $r=-DK_\Sig-1>0$.
Due to Lemma \ref{l45}, we can fix any dense subset in ${\mathcal D}_1^{\DP}$
and check the statement for the surfaces belonging to this subset.
Throughout this section, we assume that $\Sig\in U_1\cap U_2$.

We
prove the invariance of
Welschinger numbers
by studying wall-crossing events when
moving either one real point
of the given
collection, or a pair of complex conjugate points.

\smallskip

\subsubsection{Moving a real point of configuration}\label{sec-real}
Suppose that $2m<r$.
Let tuples $\bw'\cup\{w^{(0)}\},\bw'\cup\{w^{(1)}\}\in
\mathcal P_{r,m}(\Sig,F)$, where $\bw'\in
\mathcal P_{r-1,m}(\Sig,F)$, be such that
the sets ${\mathcal R}(\Sig,D,F,\bw'\cup\{w^{(0)}\})$ and
${\mathcal R}(\Sig,D,F,\bw'\cup\{w^{(1)}\})$ are finite and presented by immersions
(see Lemma \ref{ln12}). We
prove that
\begin{equation}W_m(\Sig,D,F,\varphi,\bw'\cup\{w^{(0)}\})=W_m(\Sig,D,F,\varphi,\bw'\cup\{w^{(1)}\})\
.\label{en15}\end{equation}

Due to Lemma \ref{ln16}, by a small deformation of $\bw'$ we can reach the following:
whenever an element $[\nu:\widehat C\to\Sig]\in\overline{{\mathcal M}^{br}_{0,0}}(\Sig,D)
\setminus{\mathcal M}^{br}_{0,0}(\Sig,D)$ is such that $\nu(\widehat C)\supset\bw'$,
the element $[\nu:\widehat C\to\Sig]$ satisfies the conditions of Lemma \ref{ln16}(i),
$-D_1K_\Sig-1$ points of $\bw'$ lie on $C_1\setminus(\Sing(C_1)\cup C_2)$, and
the remaining $-D_2K_\Sig-1$ points of $\bw'$ lie on $C_2\setminus(\Sing(C_2)\cup C_1)$.

There exists a smooth
real-analytic path $\sig:[0,1]\to F$ lying in the real part of some
smooth real algebraic curve $\sig(\C)\subset\Sig$,
such that $\sig$ is disjoint from all the points of $\bw'$, $\sig(0)=w^{(0)}$,
$\sig(1)=w^{(1)}$, and in the family $$\overline{{\mathcal M}_{0,r}^{br}}
(\Sig,D;\bw',\sig)=\{[\nu:\hat C\to
\Sig,\bp]\in\overline{{\mathcal M}_{0,r}^{br}}(\Sig,D)\ :\ \nu(\bp')=\bw',\
\nu(p_r)\in\sig\}\ ,$$ where $\bp'=\bp\setminus\{p_r\}$, all but
finitely
many elements
belong to ${\mathcal M}^{im}_{0,r}(\Sig,D)$, and
the remaining
elements $[\nu:\hat C\to\Sig,\bp]$ (corresponding to some values $t\in I_0\subset
[0,1]$, $|I_0|<\infty$) are
such that: \begin{enumerate} \item[(D1$_{re}$)] either $[\nu:\hat C\to\Sig]\in\overline{{\mathcal M}_{0,0}^{br}}(\Sig,D)$ is as in
Lemma \ref{ln16}(i),
the point $w^{(t)}\in \sig\cap C_2$ belongs to $C_2\setminus(\Sing(C_2)\cup C_1\cup
\bw')$,
and the germ of $\sig(\C)$ at $w^{(t)}\in C$ intersects
$C_2$ transversally; \item[(D2$_{re}$)] or
$[\nu:\hat C\to\Sig]\in{\mathcal M}_{0,0}^{br}(\Sig,D)\setminus{\mathcal M}^{im}_{0,0}(\Sig,D)$, the point $w^{(t)}\in\sig\cap C$, where $C=\nu(\hat C)$, belongs to $C\setminus(\Sing(C)\cup\bw')$, and
the germ of $\sig(\C)$ at $w^{(t)}$ intersect
$C$ transversally.
\end{enumerate}

Denote by $M_{[\nu,\bp]}$ the germ of $\overline{{\mathcal M}_{0,r}^{br}}(\Sig,D;\bw',\sig)$ at
an element $[\nu:\hat C\to\Sig,\bp]$.

If $[\nu:\hat C\to\Sig]\in{\mathcal M}^{im}_{0,0}(\Sig,D)$, or $[\nu:\hat C\to\Sig]$ satisfies condition
(D1$_{re}$), then, by Lemma \ref{p2}, the germ $M_{[\nu,\bp]}$ is diffeomorphically mapped by $\Ev$ onto the germ $(\sig,w^{(t)})$.
Moreover, the
Welschinger sign $\mu(\nu,\varphi)$ does not change along $M_{[\nu,\sig]}$. This is evident for
$[\nu:\hat C\to\Sig]\in{\mathcal M}^{im}_{0,0}(\Sig,D)$, and, under condition (D1$_{re}$), immediately follows from the fact that
$\nu$ maps the germ of $\hat C$ at the node to a pair of real smooth branches
that intersect transversally and undergo a standard smoothing in the considered bifurcation.

Under the hypotheses of condition (D2$_{re}$), the
required constancy of the Welschinger number $W_m(\Sig,D,F,\varphi,\bw'\cup\{w^{(t)}\})$
immediately follows from Lemmas \ref{leg}, \ref{new_lemma1}
and \ref{ln2}.

\subsubsection{Moving a pair of imaginary conjugate points}\label{sec-imaginary}
Assume that
$m \geq 1$.
Let tuples $\bw'\cup\{w^{(0)},\Conj\,w^{(0)}\},\bw'\cup\{w^{(1)},
\Conj\,w^{(1)}\}\in
\mathcal P_{r,m}(\Sig,F)$, where $\bw'\in
\mathcal P_{r - 2, m - 1}(\Sig,F)$,
be such that
the sets $${\mathcal R}(\Sig,D,F,\bw'\cup\{w^{(0)},\Conj\,w^{(0)}\})\quad\text{and}\quad
{\mathcal R}(\Sig,D,F,\bw'\cup\{w^{(1)},\Conj\,w^{(1)}\})$$ are finite and presented by immersions
(see Lemma \ref{ln12}). We
prove that
\begin{equation}W_m(\Sig,D,F,\varphi,\bw'\cup\{w^{(0)},\Conj\, w^{(0)}\})=W_m(\Sig,D,F,\varphi,\bw'\cup\{w^{(1)},\Conj\,w^{(1)}\})\
.\label{en16}\end{equation}

Due to Lemma \ref{ln16}, by a small deformation of $\bw'$ we can reach the following:
for any point $w$ of a certain
Zariski open subset $\Sig_{\bw'}\subset\Sig\setminus\bw'$,
whenever for an element $[\nu:\widehat C\to\Sig]\in\overline{{\mathcal M}^{br}_{0,0}}(\Sig,D)
\setminus{\mathcal M}^{br}_{0,0}(\Sig,D)$ we have $\nu(\widehat C)\supset\bw'\cup\{w\}$,
this element $[\nu:\widehat C\to\Sig]$ satisfies the conditions of Lemma \ref{ln16}(i),
$-D_1K_\Sig-1$ points of $\bw'$ lie on $C_1\setminus(\Sing(C_1)\cup C_2)$, and
the remaining $-D_2K_\Sig-2$ points of $\bw'$ and the point $w$
lie on $C_2\setminus(\Sing(C_2)\cup C_1)$.
Further on, assuming this property
of $\bw'$, we can find
a smooth
real-analytic path $\sig:[0,1]\to\Sing\setminus\R\Sig$ lying in some
smooth real algebraic curve $\sig(\C)\subset\Sig\setminus\R\Sig$,
such that $\sig$
starts
at $w^{(0)}$ and ends up at $w^{(1)}$,
avoids all the points of $\bw'$, and satisfies the following condition
({\it cf.} section \ref{sec-real}):
for all
but finitely many elements of the family
\begin{eqnarray}\overline{{\mathcal
M}_{0,r}^{br}}(\Sig,D;\bw',\sig,\Conj\,\sig)&=&\{[\nu:\hat C\to
\Sig,\bp]\in\overline{{\mathcal M}_{0,r}^{br}}(\Sig,D)\ :\nonumber\\
& &\quad\nu(\bp')=\bw',\ \nu(p_{r-1})\in\sig,\
\nu(p_r)\in\Conj\,\sig\}\ ,\nonumber\end{eqnarray} where
$\bp'=\bp\setminus\{p_{r-1},p_r\}$, we have $[\nu:\hat C\to\Sig]\in{\mathcal M}^{im}_{0,0}(\Sig,D)$,
while the
remaining elements (which correspond to some values $t\in I_0\subset[0,1]$,
$|I_0|<\infty$) are such that:
\begin{enumerate} \item[(D1$_{im}$)] either $[\nu:\hat C\to\Sig]\in\overline{{\mathcal M}_{0,0}^{br}}(\Sig,D)$ is as in
Lemma \ref{ln16}(i), where $\nu_i:\hat C_i\to\Sig$ commutes with the real structure,
$-D_1K_\Sig-1$ points of $\bw'$ lie on $C_1\setminus(\Sing(C_1)\cup C_2)$,
the remaining $-D_2K_\Sig-2$ points of $\bw'$ lie in $C_2
\setminus(\Sing(C_2)\cup C_1)$, the point $w^{(t)}\in \sig$ belongs to $C_2\setminus
(\Sing(C_2)\cup C_1\cup\bw')$,
and the germ of $\sig(\C)$ at $w^{(t)}\in C_2$ intersects
$C_2$ transversally; \item[(D2$_{im}$)] or
$[\nu:\hat C\to\Sig]\in{\mathcal M}_{0,0}^{br}(\Sig,D)\setminus{\mathcal M}^{im}_{0,0}(\Sig,D)$, the point $w^{(t)}\in\sig\cap C$, where $C=\nu(\hat C)$, belongs to $C\setminus(\Sing(C)\cup\bw')$, and
the germ of $\sig(\C)$ at $w^{(t)}$ intersects
$C$ transversally.
\end{enumerate} Notice that, in (D1$_{im}$), the case of $C_1=C_2$ is not relevant due to $-DK_\Sig>2$,
and the case of complex conjugate $C_1$
and $C_2$ does not occur either, since any real rational curve in $|D|$
must have a non-trivial one-dimensional real branch (see Section \ref{sec1.2}).

Then the proof of (\ref{en16}) literally follows the
argument of the preceding section.

\subsection{Proof of Proposition \ref{def-invariance} and Theorem \ref{t7}}
In view of Proposition \ref{p1} and Lemmas \ref{ln12}(ii) and \ref{l45}, Theorem \ref{t7} follows from Proposition \ref{def-invariance},
and,
in its turn, to prove Proposition \ref{def-invariance} it is sufficient to check the constancy of the Welschinger number in the following families:
\begin{itemize}\item
a germ of elementary deformation
$\{\Sig_t\}_{t\in(\R,0)}$, where $\Sig_0\in U_1(A_1)$,
$\Sig_t\in U_1\cap U_2$ for each $t\ne 0$, and $D_t\ne-K_{\Sig_t}$;
\item
a germ of elementary deformation
$\{\Sig_t\}_{t\in(\R,0)}$, where $\Sig_0\in{\mathcal D}_1^{\DP}\setminus U_1$,
 $\Sig_t\in U_1\cap U_2$ for each $t\ne 0$, and $D_t=-K_{\Sig_t}$.
\end{itemize}

Let $\Sig_{0}\in U_1(A_1)$,
$\Sig_t\in U_1\cap U_2$ for $t\ne 0$, and $D_t\ne-K_{\Sig_t}$.
Extend the family $\{\Sig_t\}_{t\in(\R,0)}$
to
a conjugation invariant family
$\{\Sig_t\}_{(\C,0)}$. By Lemma \ref{ln12}(2), there exists $\bw_{0}\in{\mathcal P}_{r,m}(\Sig_{0},F_{0})$
such that, for any $k\ge 0$, all elements $[\nu:\PP^1\to\Sig_{0},\bp_{0}]\in{\mathcal M}_{0,r}(\Sig_{0},D-kE,\bw_{0})$
satisfy the properties indicated in Lemma \ref{ln12}(2ii).
These elements appear only for a finite number of values of $k$ and form a finite set.
Let us associate with each of them a comb-like curve
$[\nu:\hat C\to\Sig_{0},\bp]
\in\overline{{\mathcal M}_{0,r}(\Sig_{0},D,\bw_{0})}$ such that:
\begin{itemize}\item either $\hat C\simeq\PP^1$, or $\hat C=\hat C'\cup\hat E_1\cup...\cup\hat E_k$ for some
$k>0$,
where $\hat C'\simeq\hat E_1\simeq...\simeq\hat E_k\simeq\PP^1$, $\hat E_i\cap\hat E_j=\emptyset$ for all $i\ne j$, and
$\#(\hat C'\cap \hat E_i)=1$ for all $i=1,...,k$;
\item $\bp\subset\hat C'$ and $[\nu:\hat C'\to\Sig_{0},\bp]\in{\mathcal M}^{im}_{0,r}(\Sig_{0},D-kE,\bw_{0})$, and each of $\hat E_1,...,\hat E_k$ is isomorphically
    mapped onto $E$.
\end{itemize}
Then, complement $\bw_{0}$
to a conjugation invariant family
of $r$-tuples $\bw_t\in(\Sig_t)^r$, $t\in(\C,0)$, so that $\bw_t\in{\mathcal P}_{r,m}
(\Sig_t,F_t)$ for each real $t$.
It follows from \cite[Theorem 4.2]{Va} that
each of the introduced elements $[\nu:\hat C\to\Sig_{0},\bp]
\in\overline{{\mathcal M}_{0,r}(\Sig_{0},D,\bw_{0})}$
extends
to a smooth family $[\nu_t:\hat C_t\to\Sig_t,\bp_t]\in \overline{{\mathcal M}_{0,r}^{br}}(\Sig_t,D,\bw_t)$, $t\in(\C,0)$,
where $\hat C_t\simeq\PP^1$ and $\nu_t$ is an immersion for
all $t\ne 0$, and, furthermore, each element of ${\mathcal M}_{0,r}
(\Sig_t,D,\bw_t)$, $t\in(\C,0)\setminus\{0\}$ is included into some of the above families.
Thus,
the Welschinger number $W(\Sig_t,D,F_t,\varphi_t,\bw_t)$ remains constant in $t\in(\R,0)\setminus\{0\}$,
since the only change of the topology in the real loci of the curves under
the count
consists in smoothing of non-solitary
nodes, while the difference between the homology classes of
the halves $[C_\pm(t)]$ in $H_2(\Sig_t, F_t; \Z/2)=H_2(\Sig_0, F_0;\Z/2)$ with $t<0$
and those with $t>0$ belongs to $(1+\conj_*)H_2(\Sig_0,F_0;\Z/2)$ and, hence $[C_\pm(t)]\circ\phi_t$
does not depend on $t$.

Assume that $\Sig_0\in{\mathcal D}_1^{\DP}\setminus U_1$
$\Sig_t\in U_1\cap U_2$ for $t\ne 0$,  and $D_t=-K_{\Sig_t}$.
In this case we deal with a family of real elliptic pencils $|-K_{\Sig_t}|$, $t\in(\R,0)$, such that the central one
$|-K_{\Sig_0}|$ has a real cuspidal curve $C_0\in|-K_{\Sig_{0}}|$ and, otherwise, the family is generic. As it can be seen from
the local Weierstrass normal form, due to the above genericity the  image of $|-K_{\Sig_0}|$ in
the base $(\C^2,0)$ of the versal deformation of the cuspidal point intersects the tangent space to the discriminant locus,
that is the cusp curve $27p^2+4q^3=0$ in terms of Weierstrass coordinates $p,q$,   transversally at one point.
Therefore, for $t\in (\R,0)$ on one side of $t=0$
the singular curves in $|-K_{\Sig_t}|$ close to $C_0$ form a pair of complex conjugate curves, while
for $t\in (\R,0)$ on the opposite side of $t=0$ they are real, one with a solitary node, and the other one with a cross point.
Thus, the total Welschinger number is the same on the both sides.

{\ncsc Institut de Math\'ematiques de Jussieu - Paris Rive Gauche\\[-21pt]

Universit\'e Pierre et Marie Curie\\[-21pt]

4 place Jussieu,
75252 Paris Cedex 5,
France} \\[-21pt]

{\it and} {\ncsc D\'epartement de math\'ematiques et applications\\[-21pt]

Ecole Normale Sup\'erieure\\[-21pt]

45 rue d'Ulm 75230 Paris Cedex 5, France} \\[-21pt]

{\it E-mail address}: {\ntt     ilia.itenberg@imj-prg.fr}

\vskip10pt

{\ncsc Universit\'e de Strasbourg et IRMA \\[-21pt]

7, rue Ren\'e Descartes, 67084 Strasbourg Cedex, France} \\[-21pt]

{\it E-mail address}: {\ntt kharlam@math.unistra.fr}

\vskip10pt

{\ncsc School of Mathematical Sciences \\[-21pt]

Raymond and Beverly Sackler Faculty of Exact Sciences\\[-21pt]

Tel Aviv University \\[-21pt]

Ramat Aviv, 69978 Tel Aviv, Israel} \\[-21pt]

{\it E-mail address}: {\ntt shustin@post.tau.ac.il}

\end{document}